\documentclass[11pt]{article}
\usepackage{latexsym, amsmath, amssymb}
\usepackage{amsmath}
\usepackage{amsfonts}
\usepackage{amssymb}
\usepackage{amscd}
\usepackage{amsthm}
\usepackage{fancyhdr}
\newtheorem{theorem}{Theorem}[section]
\newtheorem{proposition}[theorem]{Proposition}
\newtheorem{remark}[theorem]{Remark}

\newtheorem{definition}[theorem]{Definition}
\newtheorem{corollary}[theorem]{Corollary}
\newtheorem{example}[theorem]{Example}

\def\C{\mathbb C}

\def\Q{\mathbb Q}

\newcommand\sE{{\mathcal E}}
\newcommand\sF{{\mathcal F}}

\newcommand\sO{{\mathcal O}}

\newcommand\sS{{\mathcal S}}

  \def \tab#1{\kern #1 truein}

  \def\E{\hbox{${\cal E}$}}
  \def\F{\hbox{${\cal F}$}}
  \def\G{\hbox{${\cal G}$}}
  \def\A{\hbox{${{\cal A}}$}}
  \def\B{\hbox{${\cal B}$}}
  \def\C{\hbox{${\cal C}$}}

  \def\Q{\hbox{${\cal Q}$}}

\begin{document}
\title{A Few Splitting Criteria for Vector Bundles} 
\date{}
\author{Francesco Malaspina \\
 Dipartimento di Matematica Universit\`a di Torino\\
via Carlo Alberto 10, 10123 Torino, Italy\\ 
{\small\it e-mail: alagalois@yahoo.it, fax: $00390116702878$}}

   \maketitle
   \def\thefootnote{}
  \baselineskip=+.5cm\begin{abstract}We prove a few splitting criteria for vector bundles on a quadric hypersurface and Grassmannians. We give also some cohomological splitting conditions for rank $2$ bundles on multiprojective spaces. The tools are monads and a Beilinson's type spectral sequence generalized by Costa and Mir\'o-Roig.

      \end{abstract}

   \section*{Introduction}
A well known result by Horrocks (see \cite{Ho}) characterizes the vector bundles without intermediate cohomology on a projective space as direct sum of line bundles.
This criterion fails on  more general varieties. In fact there exist non-split vector bundles  without intermediate cohomology. This bundles are called ACM bundles.\\
On a quadric hypersurface $\Q_n$  there is a theorem that classifies all the ACM bundles (see \cite{Kn}) as  direct sums of line bundles and spinor bundles up to a twist (for generalities about spinor bundles see \cite{Ot2}).\\ 
Ottaviani has generalized Horrocks criterion to  quadrics and Grassmanniann giving cohomological splitting conditions for vector bundles (see \cite{Ot1} and \cite{Ot3}).\\
The starting point of this note is \cite{CM1} where Costa and Mir\'o-Roig give a  new proof of Horrocks and Ottaviani's criteria by using different techniques.
They also discuss that Horrocks criterion can be extended to multiprojective spaces.\\
Our aim is to use these results and techniques in order to prove  several splitting criteria for vector bundles on a quadric hypersurface and Grassmannians (in particular $G(2,n)$). All these criteria are not equivalent to those by Ottaviani. We give also some cohomological splitting conditions for rank $2$ bundles on $\mathbb{P}^n$ and multiprojective spaces. The tools are a Beilinson's type spectral sequence generalized and monads.\\
Beilinson's Theorem was stated in $1978$ and since then it has become a major tool in classifying vector bundles over projective spaces. Beilinson's spectral sequences was generalized by Kapranov to hyperquadrics and Grassmannians and by Costa and Mir\'o-Roig (see \cite{CM1}) to any smooth projective variety of dimension $n$ with an $n$-block collection.\\
A monad on $\mathbb{P}^n$ or, more generally, on a projective variety $X$, is a complex of three vector bundles
$$0 \rightarrow \A \xrightarrow{\alpha} \B \xrightarrow{\beta} \C \rightarrow 0$$
such that $\alpha$ is injective as a map of vector bundles and $\beta$ is surjective.
Monads have been studied by Horrocks, who proved (see \cite{Ho} or \cite{BH}) that every vector bundle on $\mathbb{P}^n$  is the homology of a suitable minimal monad.
This correspondence holds also on a projective variety $X$ ($\dim X\geq 3$) if we fix a very ample line bundle $\sO_X(1)$. Indeed the proof of the result in (\cite{BH} proposition $3$) can be easily extended to $X$ (see \cite{Ml})). \\
Rao, Mohan Kumar and Peterson on $\mathbb{P}^n$ (see \cite{KPR1}), and the author on quadrics have successfully used this tool to investigate the intermediate cohomology modules of a vector bundle  and give cohomological splitting conditions .\\
Here we do the same on Grassmannianns of type $G(2,n)$ and go deeper into the investigation on cohomological conditions on quadrics.\\
I wish to thank Laura Costa and Rosa-Maria Mir\'o-Roig for the warm hospitality in Barcelona and the useful discussion which gave the first hint for this paper.

\section{Preliminaries}
Throughout the paper $X$ will be a smooth projective variety defined over the complex numbers $\mathbb{C}$ and we denote by $\mathcal{D}=D^b(\sO_X-mod)$ the derived category of bounded complexes of coherent sheaves of $\sO_X$-modules.\\ 
For the notations we refer to \cite{CM1}.\\
Now we give the definition of $n$-block collection in order to introduce a Beilinson's type spectral sequence generalized:
\begin{definition}
An exceptional collection $(F_0, F_1, \dots , F_m)$ of objects of $\mathcal{D}$ (see \cite{CM1} Definition $2.1.$) is a block if $Ext^i_{\mathcal{D}}(F_j,F_k)=0$ for any $i$ and $j\not=k$.\\
An $n$-block collection of type $(\alpha_0, \alpha_1, \dots , \alpha_n)$ of objects of $\mathcal{D}$ is an exceptional collection $$(\E_0, \E_1, \dots, \E_m)=(E^0_{1},\dots, E^0_{\alpha_0}, E^1_{1},\dots, E^1_{\alpha_1}, \dots, E^n_{1},\dots, E^n_{\alpha_n})$$
such that all the subcollections $\E_i=(E^i_{1},\dots, E^i_{\alpha_i})$ are blocks.
\end{definition}
\begin{example}  $(\sO_{\mathbb{P}^n}(-n), \sO_{\mathbb{P}^n}(-n+1),\dots , \sO_{\mathbb{P}^n})$ is an $n$-block collection of type $(1, 1, \dots , 1)$ on $\mathbb{P}^n$ (see \cite{CM1} Example $2.3. (1)$).
\end{example}
\begin{example}\label{e2}  Let us consider a smooth quadric hypersurface $\Q_n$ in $\mathbb P^{n+1}$.\\
We use the unified notation $\Sigma_*$ meaning that for even $n$ both the spinor bundles $\Sigma_1$ and $\Sigma_2$ are considered, and for $n$ odd, the spinor bundle $\Sigma$.\\
$ (\E_0, \sO(-n+1),\dots , \sO(-1), \sO)$,
where $\E_0=(\Sigma_*(-n))$,  is an $n$-block collection of type $(1, 1, \dots , 1)$ if $n$ is odd, and of type $(2, 1, \dots , 1)$ if $n$ is even (see \cite{CM1} Example $3.4. (2)$).\\
Moreover we can have several $n$-block collections:
$$\sigma_j = (\sO(j),\dots ,\sO(n-1), \E_{n-j}, \sO(n+1),\dots , \sO(n+j-1))$$
where $\E_{n-j}=(\Sigma_*(n-1))$ and $1\leq j\leq n$ (see \cite{CM2} Proposition $4.4$). 
\end{example}
\begin{example}\label{e3}  Let $X=G(k,n)$ be the Grassmanniann of $k$-dimensional subspaces of the $n$-dimensional vector space. Assume $k>1$. We have the canonical exact sequence $$ 0 \to \sS\rightarrow \sO^n \rightarrow \Q \to 0,$$
where $\sS$ denote the tautological rank $k$ bundle and $\Q$ the quotient bundle.\\
Denote by $A(k,n)$ the set of locally free sheaves $\Sigma^{\alpha}\sS$ where $\alpha$ runs over diagrams fitting inside a $k\times (n-k)$ rectangle. We have $\Sigma^{(p,0, \dots ,   0)}\sS=S^p\sS$, $\Sigma^{(1, 1, 0, \dots ,0)}\sS=\wedge^2\sS=\sS(1)$ and $\Sigma^{(p_1+t, p_2+t, \dots p_k+t)}\sS=\Sigma^{(p_1, p_n, \dots, p_k)}\sS(t)$.\\
$A(k,n)$ can be totally ordered in such a way that we obtain a $k(n-k)$-block collection $\sigma=(\E_0,\dots , \E_{k(n-k)})$ by packing in the same block $\E_r$ the bundles $\Sigma^{\alpha}\sS$ with $|\alpha|=k(n-k)-r$ (see \cite{CM1} Example $3.4. (1)$).\\
\end{example}
\begin{example}\label{e3}  Let $X=\mathbf P^{n_1}\times\dots \times\mathbf P^{n_r}$ be a multiprojective space and $d=n_1+\dots+n_r$.\\
For any $0\leq j\leq d$, denote by $\E_j$ the collection of all line bundles on $X$ $$\sO_X(a_1^j, a_2^j, \dots , a_r^j)$$
with $-n_i\leq a_i^j\leq	0$ and $\sum_{i=1}^r a_i^j= j-d$.  Using the K\"unneth formula we have that $\E_j$ is a block and that $$(\E_0,\dots , \E_{d})$$ is a $d$-block collection of line bundles on $X$ (see \cite{CM1} Example $3.4. (3)$).
\end{example}
Beilinson's Theorem was generalized by Costa and Mir\'o-Roig to any smooth projective variety of dimension $n$ with an $n$-block collection of coherent sheaves which generates $\mathcal{D}$.\\
They also prove the following interesting proposition (see \cite{CM1} Proposition $4.1.$):
\begin{proposition}\label{p1}[Costa, Mir\'o-Roig] Let $X$ be a smooth projective variety of dimension $n$ with an $n$-block collection $\sigma=(\E_0,\dots , \E_{n})$, $\E_i=(E_0^i,\dots , E_{\alpha_i}^i)$ of coherent sheaves on $X$. Let assume that $\E_n=(\sO_X(e))$ where $e$ is an integer.\\
Let $\F$ be a coherent sheaf on $X$ such that for any $-n\leq p\leq -1$ and $1\leq i\leq \alpha_{n+p}$ $$H^{-p-1}(X, \F\otimes E_i^{p+n})=0.$$
Then $\F$ contains $\sO(-e)^{h^0(\sF\otimes\sO(e))}$ as a direct summand.
\end{proposition}
In the next section we will apply this proposition in order to prove splitting criteria.

\section{Splitting criteria for vector bundles}
Let $X$ be a smooth projective variety of dimension $n$ with an $n$-block collection $\sigma=(\E_0,\dots , \E_{n})$, $\E_i=(E_0^i,\dots , E_{\alpha_i}^i)$ of coherent sheaves on $X$. Let assume that $\E_n=(\sO_X(e))$ where $e$ is an integer.\\ In order to get splitting criteria for a bundle $\F$, first of all, we need $H^0(\F\otimes\sO(e))\not=0$ and $H^{0}(X, \F\otimes E_i^{n-1})=0$ for any $1\leq i\leq \alpha_{n-1}$. It's hence easier when $\E_{n-1}=(\sO(e-1))$.\\

On $\mathbf P^n$, there is a collection with $\E_n=(\sO))$ and $\E_{n-1}=(\sO(-1))$; so we can prove the following result:\\
\begin{proposition}Let $\mathcal E$ be rank two a vector bundle on $\mathbf P^n$.
Let $m$ be an integer such that $H^0(\E(m-1))=0$ and $H^0(\E(m))\not=0$.\\

Then $\mathcal E$ splits if and only if $H^i(\mathbf P^n, \E(m-1-i))=0$
for $0< i < n$.
\begin{proof} For $\E(m)$ all the conditions of (Proposition \ref{p1}) are satisfied and we can conclude that $\E(m)$ contains $\sO^{h^0(\sE(m)\otimes\sO)}$ as a direct summand. This means that $\E$ splits.
\end{proof}
\end{proposition}
\begin{remark}The same proof shows that if $\E$ is  an indecomposable vector bundle on $\mathbf P^n$ of any rank and $m$ is an integer such that $H^0(\E(m-1))=0$ but $H^0(\E(m))\not=0$.\\ Then $\E(m)\cong\sO$ if and only if $H^i(\mathbf P^n, \E(m-1-i))=0$
for $0< i < n$.\end{remark}
\begin{remark}Costa and Mir\'o-Roig re-prove Horrocks criterion by iterating this argument (see \cite{CM1} Corollary $4.2.$)
\end{remark}

On a smooth quadric hypersurface $\Q_n$ in $\mathbb P^{n+1}$ we can have many splitting criteria on $\Q_n$ (not equivalent to those by Ottaviani) by using different $n$-block collections.\\
We use the unified notation $\Sigma_*$ meaning that for even $n$ both the spinor bundles $\Sigma_1$ and $\Sigma_2$ are considered, and for $n$ odd, the spinor bundle $\Sigma$.

\begin{remark}Let $\E$ be a  vector bundle on $\Q_n$ ($n> 2$).\\
Let $m$ be an integer such that $H^0(\E(m-1))=0$ and $H^0(\E(m))\not=0$.\\
The following conditions are equivalent:
\begin{enumerate}
\item $\E$ splits into a direct sum of line bundles.
\item $H^1_*(\Q_n, \E)=\dots=H^{n-2}_*(\Q_n, \E)=H^{n-1}_*(\Q_n,
\E\otimes \Sigma_*)=0$.\\
Or, if $\E$ has rank two, $H^i(\Q_n, \E(m-1-i))=0$ for $2\leq i\leq n-2$ and\\ $H^{n-1}(\Q_n,
\E\otimes \Sigma_*(m-n))=0$.
\item There exists an integer $j$, $2\leq j\leq n-2$ such that\\ $H^1_*(\Q_n, \E)=\dots = H^{j}_*(\Q_n,
\E\otimes \Sigma_*)=\dots=H^{n-1}_*(\Q_n, \E)=0$.
\end{enumerate}
\begin{proof}$(1)\Leftrightarrow (2)$ We consider the  $n$-block collection $ (\E_0, \sO(-n+1),\dots , \sO(-1), \sO)$,
where $\E_0=(\Sigma_*(-n))$. Then we  apply (Proposition \ref{p1}) to $\E(m)$ as before and we obtain the result for rank two bundles.\\
By iterating this argument we have the general statement which is \cite{CM1} Corollary $4.3.$\\
$(1)\Leftrightarrow (3)$ We apply (Proposition \ref{p1}) by using the following $n$-block collections:
$$\sigma_j = (\sO(j),\dots ,\sO(n-1), \E_{n-j}, \sO(n+1),\dots , \sO(n+j-1))$$
where $\E_{n-j}=(\Sigma_*(n-1))$ and $3\leq j\leq n$. (see Example \ref{e2}).\\
In all these collections we have $\E_{n}=(\sO(e))$ and $\E_{n-1}=(\sO(e-1))$.\\
\end{proof}
\end{remark}
\begin{remark}If $n=2$ the only collection that we can consider is $$(\E_0=\{\sO(-2,-1)),\sO(-1,-2)\}, \sO(-1,-1), \sO).$$
So we get the following result:\\
let $\E$ be a rank $r$ bundle on $\Q_2$. Then there are $r$ integer $t_1, \dots t_r$ such that $\E\cong \bigoplus_{i=1}^r \sO(t_i,t_i)$ if and only if $H^1(\E(-2+t,-1+t)))=H^1(\E(-1+t,-2+t))=0$ for every $t\in\mathbb{Z}$.
\end{remark}
Now we introduce the following tool: the monads.\\
Let $\E$ be a vector bundle on  a nonsingular subcanonical, irreducible ACM projective variety of dimension $n$ ($n>2$). There is the corresponding minimal monad
  $$0 \rightarrow \A \xrightarrow{\alpha} \B
\xrightarrow{\beta} \C \rightarrow 0,$$ where $\A$ and $\C$ are sums of line
bundles and $\B$ satisfies:
\begin{enumerate}
\item $H^1_*(\B)=H^{n-1}_*(\B)=0$ 
\item $H^i_*(\B)=H^i_*(\E)$ \ $\forall i, 1<i<n-1
$.
\end{enumerate}
Then by (\cite{Ml} Theorem $1.2.$), if $n$ is odd and rank $\E<n-1$, or if $n$ is even and rank $\E<n$, the bundle $\B$ cannot split as a direct sum of line bundles.\\
So we can improve the above criteria in the case of bundle with a small rank:\\
\begin{theorem}Let $\E$ a vector bundle on $\Q_n$ ($n> 3$).\\

 If $n$ is odd and rank $\E<n-1$, or if $n$ is even and rank $\E<n$, then the following conditions are equivalent:
\begin{enumerate}
\item $\E$ splits into a direct sum of line bundles.
\item $H^{1}_*(\Q_n,
\E\otimes \Sigma_*)=H^2_*(\Q_n, \E)=\dots=H^{n-2}_*(\Q_n, \E)=0$.
\item $H^2_*(\Q_n, \E)=\dots=H^{n-2}_*(\Q_n, \E)=H^{n-1}_*(\Q_n,
\E\otimes \Sigma_*)=0$.
\item There exists an integer $j$, $2\leq j\leq n-2$ such that\\ $H^2_*(\Q_n, \E)=\dots = H^{j}_*(\Q_n,
\E\otimes \Sigma_*)=\dots=H^{n-2}_*(\Q_n, \E)=0$.
\end{enumerate}
\begin{proof}$(1)\Leftrightarrow (2)$ This is \cite{Ml} Theorem $2.2.$\\
$(1)\Leftrightarrow (3)$ Let us assume that $\E$ does not splits. By the above theorem we can assume that $H^1_*(\E)\not=0$ or $H^{n-1}_*(\E)\not=0$. Let us hence consider a minimal monad for $\E$

$$ 0 \to \A \xrightarrow{\alpha} \B
\xrightarrow{\beta} \C \to 0.$$  We call $\G=\ker\beta$, and 
we tensor by a spinor bundle $\Sigma_*$ the two sequences  $$ 0 \to \G \rightarrow \B
\rightarrow \C \to 0,$$ and $$ 0 \to \A \rightarrow \G \rightarrow
\E \to 0,$$ we get $$ 0 \to \G\otimes\Sigma_* \rightarrow \B\otimes\Sigma_*
\rightarrow \C\otimes\Sigma_* \to 0,$$ and $$ 0 \to \A\otimes\Sigma_*
\rightarrow \G\otimes\Sigma_* \rightarrow \E\otimes\Sigma_* \to 0.$$ So we can
see that $\forall i=1, ...,n-1$, if
$$H^i_*(\E\otimes\Sigma_*)=0,$$  we have also
$$H^i_*(\B\otimes\Sigma_*)=0.$$
For the bundle $\B$, hence all the conditions of the above theorem are
satisfied then $\B$ has to split and this is a contradiction according to (\cite{Ml} Theorem $1.2.$).\\
$(1)\Leftrightarrow (4)$ As above.
\end{proof}
\end{theorem}
\begin{remark} This means that for instance a rank two bundle $\E$ on $\Q_4$ splits if and only if $H^{2}_*(\Q_n,
\E\otimes \Sigma_*)=0$.
\end{remark}

Let $X=G(k,n)$ be the Grassmanniann of $k$-dimensional subspaces of the $n$-dimensional vector space. Assume $k>1$. 
Let us consider the $k(n-k)$-block collection $\sigma=(\E_0,\dots , \E_{k(n-k)})$  given in ( Example \ref{e3}). Here $\E_{k(n-k)}=(\sO_X)$ but $\E_{k(n-k)-1}$ is more complicated.\\
We can prove the following splitting criterion:\\
\begin{theorem}Let $\E$ be a  vector bundle on $X=G(k,n)$ and $d=k(n-k)$, then 
the following conditions are equivalent:
\begin{enumerate}
\item $\E$ splits into a direct sum of line bundles.
\item $H^i_*(X, \E\otimes\Sigma^{\alpha}\sS)=0$ when $|\alpha|=1+i$, $1\leq i\leq d-1$ and\\
$H^1_*(X, \E\otimes\Q^{\vee})=H^2_*(X, \E\otimes S^2\Q^{\vee})=\dots = H^{k-1}_*(X, \E\otimes S^{k-1}\Q^{\vee})=0$.
\end{enumerate}
\begin{proof}We may suppose that $\E$ is indecomposable.\\ Let $m$ be an integer such that $H^0(\E(m-1))=0$ and $H^0(\E(m))\not=0$.\\
In order to apply (Proposition \ref{p1}) we need $H^0(X, \E(m)\otimes\Sigma^{\alpha}\sS)=0$ when $|\alpha|=1$. Since $\Sigma^{(1,0,\dots ,0)}\sS=\sS$ we simply need $H^0(X, \E(m)\otimes\sS)=0$.\\
By using the natural isomorphism $\wedge^{k-1}\sS^{\vee}\cong\sS(1)$ we get a long exact sequence $$ 0 \rightarrow S^{k-1}\Q^{\vee}\otimes \E(m-1)
\rightarrow S^{k-2}\Q^{\vee}\otimes \sO^n\otimes \E(m-1) \rightarrow\dots \rightarrow
 \wedge^{k-1}\sO^n\otimes \E(m-1) \rightarrow\sS\otimes \E(m) \rightarrow 0.$$
 $H^1_*(X, \E\otimes\Q^{\vee})=H^2_*(X, \E\otimes S^2\Q^{\vee})=\dots = H^{k-1}_*(X, \E\otimes S^{k-1}\Q^{\vee})=0$ and $H^0(\E(m-1))=0$ imply $H^0(X, \E(m)\otimes\sS)=0$.\\
 Now we can apply (Proposition \ref{p1}) and conclude that $\E(m)\cong \sO$.
 \end{proof}
\end{theorem}
Let us study more careful the case of $k=2$.
\begin{theorem}Let $\E$ be a  vector bundle on $X=G(2,n)$ and $d=2(n-2)$, then 
the following conditions are equivalent:
\begin{enumerate}
\item $\E$ splits into a direct sum of line bundles.
\item
$H^i_*(X, \E\otimes S^{i+3-2j}\sS)=0$, $1\leq i\leq n-3$, $j>0$ and\\
$H^{d-i}_*(X, \E\otimes S^{i+1-2j}\sS)=0$, $1\leq i\leq n-2$, $j>0$.
\end{enumerate}
Moreover  if  rank $\E<d$ the conditions $H^1_*(X, \E)=0$ and $H^{n-1}_*(X, \E)=0$ are not necessary.\\
We have that the following conditions are equivalent:
\begin{enumerate}
\item $\E$ splits into a direct sum of line bundles.
\item $H^1_*(X, \E\otimes S^2\sS)=0$,\\
$H^i_*(X, \E\otimes S^{i+3-2j}\sS)=0$, $2\leq i\leq n-3$, $j>0$ and\\
$H^{d-i}_*(X, \E\otimes S^{i+1-2j}\sS)=0$, $2\leq i\leq n-2$, $j>0$.
\end{enumerate}
\begin{proof}These hypothesis are corresponding  to $H^i_*(X, \E\otimes\Sigma^{\alpha}\sS)=0$ when $|\alpha|=1+i$, $1\leq i\leq d-1$ since  $\Sigma^{(p,0)}\sS=S^p\sS$, $\Sigma^{(1,1)}\sS=\wedge^2\sS=\sS(1)$ and $\Sigma^{(p+t,q+t)}\sS=\Sigma^{(p,q)}\sS(t)$.\\
Moreover according with the above theorem we need $H^1_*(X, \E\otimes\Q^{\vee})=0$. By using the natural isomorphism $\wedge^{n-3}\Q\cong\Q^{\vee}(1)$ we get a long exact sequence $$ 0 \rightarrow S^{n-3}\sS\otimes \E
\rightarrow S^{n-4}\sS\otimes \sO^n\otimes \E \rightarrow\dots \rightarrow
 \wedge^{n-3}\sO^n\otimes \E \rightarrow\Q^{\vee}(1)\otimes \E \rightarrow 0.$$
 $H^1_*(X, \E)=H^2_*(X, \E\otimes \sS)=\dots = H^{n-2}_*(X, \E\otimes S^{n-3}\sS)=0$ (that are some of our hypothesis) imply $H^1_*(X, \E\otimes\Q^{\vee})=0$ so we don't need any additional hypothesis and we have the first part of the theorem.\\
 
Let $\E$ be now a bundle with rank $\E<d$. Let us assume that $\E$ does not splits but satisfy all our cohomological conditions. By the first part of the proof we can assume that $H^1_*(\E)\not=0$ or $H^{n-1}_*(\E)\not=0$. Let us hence consider a minimal monad for $\E$

$$ 0 \to \A \xrightarrow{\alpha} \B
\xrightarrow{\beta} \C \to 0.$$  We call $\G=\ker\beta$, and 
we tensor by a bundle $\F$ the two sequences  $$ 0 \to \G \rightarrow \B
\rightarrow \C \to 0,$$ and $$ 0 \to \A \rightarrow \G \rightarrow
\E \to 0,$$ we get $$ 0 \to \G\otimes\F \rightarrow \B\otimes\F
\rightarrow \C\otimes\F \to 0,$$ and $$ 0 \to \A\otimes\F
\rightarrow \G\otimes\F \rightarrow \E\otimes\F \to 0.$$ So we can
see that $\forall i=1, ...,n-1$, if
$$H^i_*(\E\otimes\F)=H^i_*(\A\otimes\F)=
H^i_*(\C\otimes\F)=0,$$  we have also
$$H^i_*(\B\otimes\F)=0.$$
For the bundle $\B$, hence all the conditions of the first part of the theorem are
satisfied then $\B$ has to split and this is a contradiction according to (\cite{Ml} Theorem $1.2.$).
 \end{proof}
\end{theorem}
\begin{remark} This splitting criterion is not equivalent to those by Ottaviani. If we consider $G(2,5)$ for instance the splitting conditions are the following:\\
$H^1_*(\E)=H^1_*(\E\otimes S^2 \sS)=0$\\
$H^2_*(\E\otimes S^3 \sS)=H^2_*(\E\otimes  \sS)=0$\\
 $H^3_*(\E\otimes S^2 \sS)=H^3_*(\E)=0$\\
 $H^4_*(\E\otimes  \sS)=0$\\
 $H^5_*(\E)=0$.\\
They do not imply $H^2_*(\E)=0$ appearing in the Ottaviani criterion.
\end{remark}

Let $X=\mathbf P^{n_1}\times\dots \times\mathbf P^{n_r}$ be a multiprojective space and $d=n_1+\dots+n_r$.\\
Let us consider the $d$-block collection $\sigma=(\E_0,\dots , \E_{d)})$  given in ( Example \ref{e3}). Here $\E_d=(\sO_X(0, \dots , 0))$ but $\E_{d-1}$ is more complicated.\\
We need the following definition:
\begin{definition}A bundle $\E$ on $X=\mathbf P^{n_1}\times\dots \times\mathbf P^{n_r}$ is normalized if $H^0(X, \E)\not=0$ and $H^0(X, \E\otimes \sO_X(a_1, \dots a_r))=0$ when $a_1, \dots a_r$ are non-positive integers not all vanishing.
\end{definition}
Now we can prove the following result:\\
\begin{proposition} Let $X=\mathbf P^{n_1}\times\dots \times\mathbf P^{n_r}$ be a multiprojective space and $d=n_1+\dots+n_r$.\\ Let $\E$ be  a normalized bundle on $X$. 
If, for any $i$, $1\leq i\leq d-1$, 
$$H^i(X, \E\otimes \sO(t_1,\dots, t_r))=0$$ when $-n_j\leq t_j\leq 0$ and $\sum_{j=1}^r t_j=-i-1$, then $\E$ contains $\sO_X$ as a direct summand.

\begin{proof}Let $\sigma=(\E_0,\dots , \E_d)$ be the $d$-block collection given in ( Example \ref{e3}). Our cohomological conditions  are corresponding to those of (Proposition \ref{p1}) so, since $\E_n=(\sO_X)$, we can conclude that $\E$ contains $\sO_X$ as a direct summand.
\end{proof}
\end{proposition}
\begin{corollary} Let $X=\mathbf P^{n_1}\times\dots \times\mathbf P^{n_r}$ be a multiprojective space and $d=n_1+\dots+n_r$.\\ Let $\E$ be  a rank two normalized bundle on $X$. 
If, for any $i$, $1\leq i\leq d-1$, 
$$H^i(X, \E\otimes \sO_X(t_1,\dots, t_r))=0$$ when $-n_j\leq t_j\leq 0$ and $\sum_{j=1}^r t_j=-i-1$, then $\E$ splits as a direct sum of line bundles.
\end{corollary}
\bibliographystyle{amsplain}

\begin{thebibliography}{99}
\bibitem
{BH}{\sc W. Barth, K. Hulek},
\emph{Monads and moduli of vector bundles}, 1978, Manuscripta Math. 25, 323-447.
\bibitem
{Be} {\sc A. A. Beilinson},
\emph{Coherent sheaves on $\mathbb{P}^n$ and Problems of Linear Algebra}, 1979, Funkt. Anal. Appl. 12, 214-216.
\bibitem
{CM1} {\sc L. Costa and R.M. Mir\'o-Roig},
\emph{Cohomological characterization of vector bundles on multiprojective spaces}, 2005, J. of Algebra, 294,    73-96.
\bibitem
{CM2} {\sc L. Costa and R.M. Mir\'o-Roig},
\emph{Monads and regularity of vector bundles on projective varieties}, 2007, Preprint. 
\bibitem
{Ho} {\sc G. Horrocks},
\emph{Vector bundles on the punctured spectrum of a ring}, 1964, Proc. London Math. Soc. (3) 14, 689-713.
\bibitem
{Ka1} {\sc M. M. Kapranov},
\emph{On the derived category of coerent sheaves on Grassmann manifolds}, 1985, Math. USSR Izvestiya, 24, 183-192.
\bibitem
{Ka2} {\sc M. M. Kapranov},
\emph{On the derived category of coerent sheaves on some homogeneous spaces}, 1988, Invent. Math., 92, 479-508.
\bibitem
{Kn} {\sc H. Kn\"{o}rrer},
\emph{Cohen-Macaulay modules of hypersurface singularities I}, 1987, Invent. Math. 88, 153-164.
\bibitem
{KPR1}{\sc N.Mohan Kumar, C. Peterson and A.P. Rao},
\emph{Monads on projective spaces}, 2003, Manuscripta Math. 112, 183-189.
\bibitem
{Ml}{\sc F. Malaspina},
\emph{Monads and Vector Bundles on Quadrics}, 2006, Preprint math.AG/0612512.
\bibitem
{Ot1}{\sc G. Ottaviani}, 
\emph{Criteres de scindage pour les fibres vectoriel sur les grassmanniennes et les quadriques }, 1987, C. R. Acad. Sci. Paris, 305, 257-260.
\bibitem
{Ot2}{\sc G. Ottaviani},
\emph{Spinor bundles on Quadrics}, 1988, Trans. Am. Math. Soc:, 307, no 1, 301-316.
\bibitem
{Ot3}{\sc G. Ottaviani},
\emph{ Some extension of Horrocks criterion to vector bundles on Grassmannians and quadrics}, 1989, Annali Mat. Pura Appl. (IV) 155, 317-341.
\end{thebibliography}

\end{document}